\documentclass[article]{amsart}
\usepackage{amsfonts}
\usepackage{graphicx}
\usepackage{amscd}
\usepackage{amsmath}
\theoremstyle{plain} 

\newtheorem{definition}{\bf Definition}

\newtheorem{lemma}{\bf Lemma}

\newtheorem{theorem}{\bf Theorem}[section]
\numberwithin{equation}{section}

\newcommand{\grad}{\operatorname{grad}}
\newcommand{\diver}{\operatorname{div}}

\newcommand{\Hess}{{\rm Hess}\,}
\newcommand{\Div}{{\rm div}\,}

\newcommand{\flux}{{\rm flux}\,}

\begin{document}
\title[Stochastic Properties of the Laplacian on Riemannian Submersions]{Stochastic Properties of the Laplacian on Riemannian Submersions  }
\author{M. Cristiane  Brand\~ao, Jobson Q. Oliveira }
\address{Departamento de Matem\'atica-UFC\\
60455-760-Fortaleza-CE-Br} \email{jobsonqo@gmail.com}\email{crismbrandao@yahoo.com.br}
\urladdr{http://www.mat.ufc.br/}

%\subjclass{Primary 53D10; Secondary 53D35, 53C40}


\begin{abstract} Based on ideas of Pigolla and Setti \cite{PS} we prove that immersed submanifolds
with bounded mean curvature of Cartan-Hadamard  manifolds are
Feller. We also  consider Riemannian submersions $\pi \colon  M \to N$ with compact
minimal fibers, and based on various criteria for parabolicity and stochastic completeness, see \cite{Grygor'yan},  we prove that $M$ is   Feller,   parabolic or stochastically complete if and only if the base $N$ is Feller, parabolic or stochastically complete  respectively.
\end{abstract}
\maketitle
\section{\bf Introduction}

Let $M$ be a  geodesically complete  Riemannian
manifold and $\triangle=\diver \circ \grad $  the Laplace-Beltrami
operator acting on  the space $C_{0}^{\infty}(M)$ of smooth
functions with compact support. The  operator $\triangle$ is
symmetric with respect to the $L^{2}(M)$-scalar product and it has a
unique  self-adjoint extension
 to a semi-bounded operator, also denoted by $\triangle$, whose domain is the set $W^{2}_{0}(M)=\{ f\in W_{0}^{1}(M),\, \triangle\! f\in L^{2}(M)\}$, see details in  \cite{davies}, where $W_{0}^{1}(M)$ is the closure  of $C_{0}^{\infty}(M)$ with respect to the norm $$ (u,v)_{1}=\int_{M}u\, v \,d\mu + \int_{M}\langle \nabla u, \, \nabla v \rangle \,d\mu.$$
 The  operator $\triangle$ defines the heat semi-group $\{e^{t\triangle}\}_{t\geq 0}$,  a family of positive definite bounded self-adjoint operators in $L^{2}(M)$ such that for any $u_{0}\in L^{2}(M)$, the function $u(x,t)\colon\!\!=(e^{t\triangle }u_{0})(x)\in C^{\infty}((0, \infty)\times M)$ solves the heat equation \begin{equation}\left\{\begin{array}{rcl}\displaystyle\frac{\partial }{\partial \,t}u(t,x)&=&\triangle_{x}u(t,x)\\
 &&\\
 u(t,x)&\stackrel{L^{2}(M)}{\longrightarrow }&u_{0}(x)\,\,\;\;\;{\rm as}\,\,\;\;\; t\to 0^{+}\end{array}\right.
\end{equation} Moreover, there exists a smooth function $p\in C^{\infty}(\mathbb{R}_{+}\times M \times M)$, called the heat kernel of $M$, such that \begin{equation} e^{t\triangle}u(x)=\int_{M}p(t,x,y)\,u(y)\,d\mu_{y},\end{equation}see \cite{dodziuk}, \cite{G-book}. In   \cite{azencott},
Azencott  studied (among other things) Riemannian manifolds  such  that the heat
semi-group  $e^{t\triangle }$ preserves the set of continuous
function vanishing at infinity. He introduced the concept
of Feller manifolds in the following definition. \begin{definition}A
complete Riemannian manifold $M$ is Feller or enjoys the Feller property for the Laplace-Beltrami operator  if
\begin{equation} e^{t\triangle }(C_{0}(M))\subset
C_{0}(M)\end{equation} where $C_{0}(M)=\{u\colon M\to \mathbb{R},
\,{\rm continous}\colon u(x)\to 0 \, \,{\rm as}\,\, x\to \infty\}.$
\end{definition}

 On the other hand, it is well known that the heat kernel has the following properties,
\begin{equation}\begin{array}{rll}\displaystyle \frac{\partial p}{\partial t} &=& \triangle_{y}p. \\
&&\\
p(t, x,y) &> &0,\\
 &&\\
 p(t, x,y)&=&\int_{M}p(s, x,z)p(t-s,z,y)dz,\\&& \\ \int_{M}p(t, x,y)dy &\leq &1\end{array}\end{equation} for all $x\in M$ and all $t>0$, $s\in (0, t)$.
From these  properties one can construct a Markov process $X_{t}$  on $M$, called  Brownian motion  on $M$, with transition density $p(x,y,t)$, see \cite[p.143]{Grygor'yan}. The corresponding measure in the space of all paths issuing from a point $x$ is denoted by $\mathbb{P}_{x}$. If  $X_{0}=x$ and  $U\subset M$ is an open set. Then $$\mathbb{P}_{x}\left( \{ X_{t}\in U\}\right)=\int_{U}p(x, y, t)dy$$ The process $X_{t}$ is stochastically complete it the total probability of the particle being found in $M$  is equal to $1$. This motivates the following definition\begin{definition}A Riemannian manifold $M$ is said to be stochastically complete if  for some, equivalently for all, $(x,t)\in M\times (0, \infty)$, $$ \int_{M}p(x, y, t)dy =1.$$
\end{definition}

In this probabilistic
point of view, Azencott \cite{azencott} proved the following
probabilistic characterization of Feller manifolds.
\begin{theorem}[Azencott]A Riemannian manifold $M$ is Feller  if and
only if, for every compact $K\subset M$ and for every $t_{0}>0$, the
Brownian motion $X_{t}\in M$ starting at $X_{0}=x_{0}$ enters in $K$
before time $t<t_{0}$ with probability that tends to zero as
$x_{0}\to \infty$.
\end{theorem}

After Azencott's paper, many authors \cite{dodziuk},  \cite{hsu-1},
\cite{hsu-2}, \cite{LK}, \cite{pinsky}, \cite{yau-1} contributed to
the theory of Feller manifolds setting geometric conditions implying
the Feller property. Most of  those geometric conditions are all on
the Ricci curvature of the manifolds although the methods employed
differs, ranging from parabolic equations to probability methods. An interesting  approach  was taken recently by Pigola and Setti in \cite{PS}, they used a characterization of minimal
solutions of certain elliptic problems due to Azencott
\cite{azencott} to set up a very useful  criteria to prove the Feller property of Riemannian manifolds(Theorem \ref{supersol compar}), similar to
those used to prove parabolicity and stochastic completeness of
Riemannian manifolds.

Stemming from \cite{PS}, the paper \cite{BPS}  considers Riemannian  manifolds that are stochastically complete and Feller simultaneously and studies solutions of certain PDE's out of a compact set and prove a number of geometric applications, see
\cite[Thms. 16 18, 20, 21 ]{BPS}.
There are many examples  of manifolds that are   Feller and stochastically complete, like the
 Cartan-Hadamard manifolds with sectional curvature with  quadratic decay, the Ricci solitons,  the properly immersed minimal submanifolds of Cartan-Hadamard manifolds, etc.
In order to apply the machinery developed in \cite{BPS}, it is  important to establish geometric criteria to ensure stochastic completeness and Feller property of Riemannian manifolds.
 In our first result we show that  any properly immersed  submanifolds of a Hadamard-Cartan manifold with bounded mean curvature vector is Feller. It is known that properly immersed  submanifolds with bounded mean curvature vector  are stochastically complete, \cite{PRS}. We also prove stochastic completeness and Feller property of  an important class of Riemannian manifolds, the Riemannian submersions  with compact minimal fibers. The Riemannian submersions were introduced by O'Neill  \cite{One66}, \cite{One67} and A. Gray \cite{Gra67} in order to produce new examples of non-negative sectional curvature manifolds, positive Ricci curvature manifolds, as a laboratory to test conjectures. Examples of  Riemannian submersions  are the coverings spaces $\pi\colon \widetilde{M}\to M$, warped product manifolds $\pi\colon (X\times_{\psi}Y, dX^{2}+\psi^{2}(x,y)dY^{2})\to X$. To give examples of Riemannian submersions with minimal fibers, let $G$ be Lie group  endowed with a bi-invariant metric and  $K$ be a closed subgroup, then the natural projection  $\pi\colon  G\to G/K$ is a Riemannian submersion with totally geodesic fibers diffeomorphic to $K$. Other examples  are the homogeneous $3$-dimensional
Riemannian manifolds with isometry group of dimension four described in details in \cite{Sco}.

   In our second result, we show that if $\pi\colon M \to N$ is a Riemannian submersion with compact minimal fibers $F$, then  $M$ is, respectively Feller, stochastically complete or parabolic if and only if $N$ is Feller, stochastically complete or parabolic.

\section{\bf Statement of the results}
  Pigola and Setti  \cite{PS}, as
consequence of the relations between the Faber-Krahn isoperimetric
inequalities and Feller property, proved the following result.
\begin{theorem}[Pigola-Setti]  Let  $\varphi\colon M
\hookrightarrow N$ be an immersion of a $m$-submanifold with mean
curvature vector $H$ into  a Cartan-Hadamard $n$-manifold $N$. If
\begin{equation} \int_{M}\vert H\vert^{m}d\mu_{M} < \infty\end{equation} then $M$ is Feller. In
particular
\begin{itemize}
\item[a.] Every Cartan-Hadamard manifold is Feller.
\item[b.]Every complete minimal submanifold of  a Cartan-Hadamard
manifold is Feller.
\end{itemize}\label{thmPS}
\end{theorem}
In our  first result we substitute the condition $ \Vert
H\Vert_{L^{m}(M)} < \infty$ in Theorem \ref{thmPS} by properness of
the immersion and boundedness of the mean curvature vector. We prove
the following theorem.

\begin{theorem} \label{Theorem 2.1}Let  $\varphi\colon M
\hookrightarrow N$ be an proper immersion of a $m$-submanifold with mean
curvature vector $H$ into  a Cartan-Hadamard $n$-manifold $N$. If $\varphi $ has bounded  mean curvature vector,  $\sup_{M}\vert H\vert <\infty$,  then $M$ is Feller.
\end{theorem} We should remark that properly immersed submanifolds of Cartan-Hadamard manifolds with  mean curvature vector with controlled growth\footnote{Meaning that $\sup_{B_{N}(p,t)\cap \varphi (M)} \vert H\vert \leq c^{2}\cdot t^{2}\cdot \log^{2}(t+2)$, $c$ constant and $t>>  1.$} are stochastically complete, see details in \cite{PRS}.

To put our second result in context let us consider a Riemannian covering $\pi\colon \widetilde{M}\to M$. It is known that $\widetilde{M}$ is stochastically complete if and only if $M$ is stochastically complete. A proof of that based on the fact that Brownian paths in $M$ lifts to Brownian paths in $\widetilde{M}$ and Brownian paths in $\widetilde{M}$ projects into Brownian paths in $M$ can be found in Elworthy's book \cite{elworthy}.  For parabolicity, the situation is different. If $\widetilde{M}$ is parabolic then $M$ is parabolic however the converse is not true in general, as observed in \cite[p.24]{PS}, the double punctured disc is parabolic and it is covered by the Poincar\`{e} disc which is not parabolic.

In our next theorem we consider parabolicity and stochastic completeness  on Riemannian submersions  $\pi \colon M\to N$ with minimal fibers $F_{p}=\pi^{-1}(p)$, $p\in N$.
\begin{theorem}\label{Theorem 2.2} Let $\pi \colon M \to  N$ be a Riemannian submersion
with minimal fibers $F_{p}=\pi^{-1}(p)$, $p\in N$. Then
\begin{itemize}\item[i.] If $M$ is parabolic then $N$ is parabolic. 
\item[ii.] If $M$ is stochastically complete then $N$ is stochastically complete.
\end{itemize}If in  addition to minimality,  the fibers $F_{p}$ are compact then we have.
\begin{itemize}\item[iii.]If $N$ is parabolic then $M$ is parabolic.
\item[iv.] If $N$ is stochastically complete then $M$ is stochastically complete.
\end{itemize}
\end{theorem}

\noindent \textbf{Observations.}
\begin{itemize}\item A Riemannian covering is a particular example of a Riemannian submersion with minimal fibers, thus the items i. and ii. extend the well known facts about parabolicity and stochastic completeness cited above.
\item The compactness of the fibers in items iii. and iv. can not be removed as one can see in the following examples.
\item[] \begin{itemize}\item[1.]$\pi\colon \mathbb{R}^{3}\to \mathbb{R}^{2}$ is a Riemannian submersion with non-compact minimal fibers $\mathbb{R}$. The base $\mathbb{R}^{2}$ is parabolic while $\mathbb{R}^{3}$ is not.
    \item[2.]  Let $M_1$, $M_{2}$ be   stochastically incomplete and stochastically complete Riemannian manifolds respectively. The projection
     $\pi\colon M_{1}\times M_{2}\to M_{2}$ is a Riemannian submersion with totally geodesic fibers $F\approx M_{1}$. The base space $M_{2}$ is stochastically complete while the total space $M_{1}\times M_{2}$ is not.
\end{itemize}
\end{itemize} Regarding the Feller property, Pigola and Setti proved the following result.
\begin{theorem}[Pigola-Setti] Let $\pi\colon \widetilde{M}\to M$ be a $k$-folding Riemannian covering, $k<\infty$. Then $\widetilde{M}$ is Feller if and only if $M$ is Feller. \label{thmFellerCov}
\end{theorem} Moreover, they show an example of an $\infty$-covering $\pi\colon \widetilde{M}\to M$ such that $\widetilde{M}$ is Feller while $M$ is not. However, they prove that if $M$ is Feller then any $k$-folding Riemannian covering, $k\leq \infty$  $\widetilde{M}$ is Feller, see
\cite[thm. 9.5]{PS}.
Our third result is an extension of Pigola-Setti's Theorem \eqref{thmFellerCov}, however it does not extend Theorem (9.5) of \cite{PS}.  We prove the following  theorem.
\begin{theorem}\label{thmFellerSub} Let $\pi \colon M \to  N$ be a Riemannian submersion with compact minimal fibers $F$.
 Then $M$ is Feller if and only if $N$ is Feller.
\end{theorem}

\section{\bf Proof of the Results}

 Let $ \varphi : M \hookrightarrow N$ an isometric immersion of a Riemannian $m$-manifold $M$ into a Riemannian $n$-manifold $N$. Let $g: N \rightarrow \mathbb{R}$ be a smooth function and  consider the function $f = g \circ \varphi$. It is well known that, (identifying $X$ with $d\varphi X$),
\begin{equation*}
  \Hess f(p)(X,Y) = \Hess g(\varphi (p))(X,Y)
+ \langle \alpha (X,Y),
\grad g \rangle (\varphi (p)),\,\, \forall \,X,Y\in T_{p}M \label{equation 2} \\
\end{equation*} Taking an orthonormal basis $\{X_{1}, \ldots, X_{m}\}$ of $T_{p}M$ and taking the trace we obtain
 \begin{equation}
\triangle f(p) = \sum_{i=1}^{m}\Hess g(\varphi(x))(X_{i},X_{i}) +
\langle H, \grad g \rangle(\varphi(p)) \label{equation 3}
\end{equation}
where $H={\rm Trace}\alpha$ denotes the mean curvature vector, see \cite{JK}.
\subsection{Proof of Theorem \ref{Theorem 2.1}}The theorem below is due to Azencott \cite{azencott}, see \cite{PS}. It relates
the Feller property and the decay at infinity
of a minimal solution of a certain Dirichlet problem.
\begin{theorem}[Azencott]\label{feller equiv} The following statements are equivalents.
\begin{enumerate}
\item[a.] M is Feller.
\item[b.] For any $\Omega \subset \subset M$ with smooth
boundary and for any constant $\lambda > 0$, the
minimal solution $h: M \setminus \Omega \rightarrow \mathbb{R}$ of
the problem

\begin{equation}\label{eqFeller}
\left\{
\begin{array}{lll}
\Delta h  =  \lambda h, & \mbox{on} & M \setminus \Omega \\
h  =  1, & \mbox{on} & \partial \Omega \\
h  >  0,  & \mbox{on} & M \setminus \Omega \\
\end{array}
\right.
\end{equation}
\end{enumerate}
Satisfies $h(x) \rightarrow 0$, as $x \rightarrow \infty$
\end{theorem}
 \noindent The  minimal positive solution $h$ for the problem \eqref{eqFeller} always exists, see \cite{azencott}. \begin{definition}We say that $u\colon M\setminus \Omega\to \mathbb{R}$ is a super-solution of the exterior Dirichlet problem \eqref{eqFeller}
if $u$ satisfies
\begin{equation}\label{eqFeller2}
\left\{
\begin{array}{rllll}
\Delta u & \leq & \lambda u, & \mbox{on} & M \setminus \Omega \\
u & \geq  &1, & \mbox{on} & \partial \Omega \\
\end{array}
\right.
\end{equation}A sub-solution is similarly defined, reversing the
inequalities in \eqref{eqFeller2}.
\end{definition}

This next theorem due to Pigola and Setti \cite{PS} establish a comparison between the
solution and the super-solution of the Dirichlet problem \eqref{eqFeller}.
\begin{theorem}[Pigola-Setti]\label{supersol compar} Let $\Omega$ a relatively compact open set with smooth
boundary $\partial \Omega$ in a Riemannian manifold $M$ and let
$ \lambda > 0$. Let $u$ and $h$ be  a positive
super-solution and a minimal solution of the problem\eqref{eqFeller} respectively.
 Then
$$h(x) \leq u(x), \,\, \forall x \in  M \setminus \Omega .$$ \\
In particular if $u(x) \rightarrow 0$ as $x \rightarrow \infty$ then
$M$ is Feller.
\end{theorem} Using Theorem \eqref{supersol compar} we prove Theorem \eqref{Theorem 2.1}.
\begin{proof} Let $p\in \varphi (M)\subset N$ and let  $\rho_{N}(x)={\rm dist}_{N}(p,x)$  be the distance function in $N$.   Let $\lambda, R > 0$  be positive constants and define $G\colon N\setminus B_{N}(p,R)\to \mathbb{R}$ given by $G(x)=g(\rho_{N})(x)$, where
$g: [R, +\infty) \rightarrow \mathbb{R}$ is given by $$g(t)=
\it{e}^{-\sqrt{\lambda}(t-R)}$$ and $B_{N}(p,R)$ is the geodesic ball of radius $R$ and center at $p$. Let $\Omega =\varphi (B_{N}(p,R))$ be a relatively compact open subset of $M$, (recall that $\varphi$ is a proper immersion) and define $u\colon M\setminus \Omega\to \mathbb{R}$ given by $u=G\circ \varphi$.  Let $x\in M$ such that $\varphi (x)\in N\setminus B_{N}(p,R)$. By the Formula \ref{equation 3} we have, taking a orthonormal basis
for $T_{\varphi(x)}M$ we have
%\begin{equation}
\begin{eqnarray*}
\triangle u(x) & = & \sum_{i=1}^{m}\Hess (g \circ
\rho_{N})(\varphi(x))(e_{i},e_{i}) + \langle H, \grad
(g \circ \rho_{N}) \rangle(\varphi(x)) \\
           % & = & \sum_{i=1}^{k}\Hess(h \circ \rho)(\varphi(x))(e_{i},
           % e_{i}) + \sum_{i=k+1}^{m}\Hess(h \circ \rho)(\varphi(x))(e_{i},
           % e_{i}) + \langle m\overrightarrow{H}, \grad
%(h \circ \rho) \rangle|_{\varphi(x)} \\
\end{eqnarray*}
%\end{equation}
Let $t=\rho_{N}(\varphi (x))$ and
choosing the orthonormal basis $\{ e_{i} \}$ for $T_{\varphi (x)}M$ such that $e_{2},
\ldots, e_{m}$ are tangent to the sphere $\partial B_{N}(p,t)$ and
$e_{1} = a\cdot (\partial/\partial t)  +
b\cdot (\partial/ \partial\theta)$, $a^{2}+b^{2}=1$,  where
$\partial/ \partial\theta \in [[e_{2}, \ldots, e_{m}]]$,
$\vert \partial/\partial \theta\vert =1$, $ \partial/\partial t =
\grad \rho_{N}$  we obtain
\begin{eqnarray}\label{eqFeller3}
\triangle u(x)& = & \sum_{i=1}^{m}\Hess(g \circ
\rho_{N})(\varphi(x))(e_{i},
            e_{i}) + \langle H, \grad
(g \circ \rho) \rangle(\varphi(x)) \nonumber \\
           & = & a^{2}g''(t) + b^{2}g'(t)\Hess \rho_{N}(\varphi(x))(\frac{\partial}{\partial\theta},\frac{\partial}{\partial\theta}
           ) \nonumber\\
           && + g'(t)\sum_{i=2}^{m}\Hess \rho_{N}(\varphi(x))(e_{i},
           e_{i}) +  g'(t) \langle \grad \rho_{N}, H
           \rangle \nonumber\\
           & = & (1-b^{2})g''(t) + b^{2}g'(t)\Hess \rho_{N}(\varphi(x))(\frac{\partial}{\partial\theta},\frac{\partial}{\partial\theta}
           )\nonumber\\ && + g'(t)\sum_{i=2}^{m}\Hess \rho_{N}(\varphi(x))(e_{i},
           e_{i}) +  g'(t) \langle \grad \rho_{N}, H
           \rangle \nonumber\\
           & \leq & g''(t) + g'(t)\langle \grad \rho_{N}, H
           \rangle
\end{eqnarray}

Observe that  $g> 0$ is positive,
$g'=
-\sqrt{\lambda}\,g<0$ and $ g'' = \lambda\, g>0.$  Thus

\begin{eqnarray*}
\triangle u(x) & \leq &  g''(t) + g'(t) \langle \grad \rho,
           H\rangle  \\
           & = & \lambda g(t) + (-\sqrt{\lambda})g(t)\langle \grad \rho,
           H\rangle  \\
           & \leq & (\lambda + \sqrt{\lambda}\displaystyle\sup_{M}\cdot \vert H\vert )g(t) \\
           & = & \mu \cdot g(t) \\
           & = & \mu \cdot u(x)
\end{eqnarray*} Moreover, if $x\in \partial \Omega$ then $u(x)=1$ and when $x \rightarrow \infty$ in $M$ then $\varphi (x)\to \infty $ in $N$. Therefore $u(x)\to 0$ as $x\to \infty$.
Let $h > 0$ the minimal solution of the problem
 $$
\left\{
\begin{array}{rllll}
\Delta h  & =&   \mu\cdot  h, & \mbox{on} & M \setminus \Omega \\
h  & =&  1, & \mbox{on} & \partial \Omega \\
h & >&  0,  & \mbox{on} & M \setminus \Omega \\
\end{array}
\right.
$$
Accordingly to Theorem \ref{supersol compar}  $$h(x) \leq u(x), \,\, \forall x \in M \setminus
 \Omega$$
 Taking, in the inequality above, the limit when $x \rightarrow \infty$ we obtain
 $$0 \leq \displaystyle\lim_{x \to \infty}h(x) \leq \displaystyle\lim_{x \to \infty}u(x)=0$$
and we conclude that $M$ is
 Feller.
\end{proof}
%\section{\bf Riemannian submersions}
\subsubsection{\bf Riemannian Submersions}\label{sub:notterm}In this section we discuss basic facts related to Riemannian submersions needed in the proof of our results, see more details in \cite{One66}. Let
 $M$ and $N$ be Riemannian manifolds, a smooth surjective map $\pi\colon M\to N$   is a {\em submersion}
if the differential $\mathrm d\pi (q)$ has maximal rank for every $q\in M$.
If $\pi:M\to N$ is a submersion, then for all
$p\in N$ the inverse image $F_p=\pi^{-1}(p)$ is a smooth embedded submanifold of $M$, that
called the \emph{fiber} at $p$.

\begin{definition}A submersion $\pi:M\to N$ is called a \emph{Riemannian submersion} if for all $p\in N$ and
all $q\in F_p$, the restriction of $\mathrm d\pi(q)$ to the orthogonal subspace
$T_q  F_p^\perp$ is an isometry onto $T_pM$.
\end{definition}

Given $p\in N$ and $q\in F_p$, a tangent vector $\xi\in T_qM$ is said to be \emph{vertical} if it
is tangent to $ F_p$ and it is said to be \emph{horizontal} if it belongs to the orthogonal space
$(T_q F_p)^\perp$.
Given $\xi\in TM$, its horizontal and vertical components are denoted respectively
by $\xi^{h}$ and $\xi^{v}$. The second fundamental form of the fibers is a symmetric tensor $\mathcal S^F:\mathcal D^\perp\times\mathcal D^\perp\to\mathcal D$, defined by
\[\mathcal S^F(v,w)=(\displaystyle\nabla^{M}_vW)^{h},\]
where  $W$ is a vertical extension of $w$ and $\nabla^M$ is the Levi--Civita connection of $M$.

For any given  vector field $X\in\mathfrak X(N)$,  there exists a unique horizontal vector field
$\widetilde{X}\in\mathfrak X(M)$  which is $\pi$-related
to $X$, this is,  for any $p\in N$ and $q\in F_p$,
then $\mathrm{ d}\pi_{q}(\widetilde{X}_{q})=X_{p}$, called {\em horizontal lifting} of $X$.  On the other hand, a horizontal vector field $\widetilde{X}\in\mathfrak X(M)$ is called \emph{basic} if it is $\pi$-related to
some vector field $X\in\mathfrak X(N)$.

Observe that the fibers are totally geodesic
submanifolds of $M$ exactly when $\mathcal S^F=0$.
The \emph{mean curvature} vector of the fiber is the horizontal vector field $H$ defined\footnote{Sometimes the mean curvature vector is defined as $H(q) = \sum_{i=1}^{k}\mathcal S^F(q)(e_{i}, e_{i})$} by
\begin{equation}H(q) =-\sum_{i=1}^{k}\mathcal S^F(q)(e_{i}, e_{i})=-\sum_{i=1}^{k}(\displaystyle\nabla^{M}_{e_{i}}e_{i})^{h}\label{eq:defmeancurvature}
\end{equation}
where  $(e_i)_{i=1}^k$ is a local orthonormal frame for the fiber
through $q$.  Observe that $H$ is not basic in general. For instance, when  the fibers
are hypersurfaces of $M$, then $H$ is basic if and only if all the fibers have constant
mean curvature. The fibers are \emph{minimal} submanifolds of $M$ when $H\equiv0$.
 The following lemma, whose proof can be found in \cite{BMP},  will play an
important role in the proof of Theorem \ref{Theorem 2.1}.

\begin{lemma}[Main]\label{Lemma 1} Let $\widetilde{X} \in \mathfrak{X}(M)$ be
a basic vector field, $\pi$-related to $X \in \mathfrak{X}(N)$. Then the
following relation between the divergence of
$\widetilde{X}$ and of $X$ at $x \in N$ and  at $\widetilde{x} \in
\mathcal{F}_{x}$ respectively,  holds.
  \begin{eqnarray}
   \Div_{M}(\widetilde{X})(\widetilde{x}) & = &
\Div_{N}(X)(x) + g_{N}(\widetilde{X}_{\widetilde{x}},
H_{\widetilde{x}}) \nonumber \\
& = & \Div_{N}(X)(x) + g_{N}(d
\pi_{\widetilde{x}}(\widetilde{X}_{\widetilde{x}})
,d \pi_{\widetilde{x}}(H_{\widetilde{x}})) \nonumber
     \end{eqnarray}
  If the fibers are minimal, then $\Div_{M}\widetilde{X} = \Div_{N}X$
\end{lemma}
 Let $u\colon N\to \mathbb{R}$ be a smooth function and denote by $\widetilde{u}=u \circ \pi \colon M\to \mathbb{R}$ its lifting to $M$. It is easy to show that $\widetilde{\grad_{N} u}= \grad_{M}\widetilde{u}$, the horizontal lifting of $\grad_{N}u$ is the gradient   of the horizontal lifting $\widetilde{u}$, $\grad_{M}\widetilde{u}$. We are denoting with a tilde superscript  $\widetilde{X}, \widetilde{u}$ the horizontal lifting of $X$, $u$, respectively.
\subsection{\bf Proof of Theorem \ref{Theorem 2.2}, items i. and ii.} The proof of items i. and ii. follows from a characterization of parabolicity and stochastic completeness in terms of the weak Omori-Yau maximum principle at infinity, proved by Pigola-Rigoli-Setti in \cite{PRS-0}, \cite{PRS}. Precisely they proved the following theorem.
\begin{theorem}[Pigola-Rigoli-Setti]\label{thm-PRS-parabolicity} A  Riemannian manifold $M$ is parabolic, $($resp. stochastically complete$)$, if and only if for every $u\in C^{2}(M)$, $u^{\ast}=\sup_{M}u<\infty$, and for every $\eta>0$ one has \begin{equation}\label{eqPRS-parabolicity}
\inf_{\Omega_{\eta}}\triangle u <0,\,\,\, (resp.\,\leq  0)
\end{equation}where $\Omega_{\eta}=\{u>u^{\ast}-\eta\}$.
\end{theorem}
Let $\pi\colon M\to N$ be a Riemannian submersion with minimal fibers $F$ where  $M$ is parabolic (resp. stochastically complete). Let us suppose, by contradiction, that $N$ is non-parabolic (resp. stochastically incomplete). By Theorem (3.3) there exists $\eta>0$ and  a  $u\in C^{2}(N)$ with $u^{\ast}<\infty$ such that $\inf_{\Omega_{\eta}}\triangle u \geq  0,\,\,\, (resp.\,> 0)$.

Let $\widetilde{u}\in C^{2}(M)$ be the horizontal lifting of $u$. Applying Lemma \ref{Lemma 1} to $X=\grad u$ and $\widetilde{X}=\grad \widetilde{u}$ one has  that $\diver_{N}X=\triangle_{N} u (x)=\diver_{M}\widetilde{X}=\triangle_{M}\widetilde{u}(y)$ for any $y\in F_{x}=\pi^{-1}(x)$.

It is cleat that $\widetilde{u}^{\ast}=\sup_{M}\widetilde{u}=u^{\ast}<\infty$ and defining $\widetilde{\Omega}_{\eta}=\{\widetilde{u}>\widetilde{u}^{\ast}-\eta\}$ one has that $$\inf_{\widetilde{\Omega}_{\eta}}\triangle_{M}\widetilde{u}= \inf_{\Omega_{\eta}}\triangle u \geq 0,\,\, (resp. \,>0) $$ showing that $M$ is non-parabolic, (resp. stochastically incomplete) contradicting the hypothesis that $M$ is parabolic, (resp. stochastically complete).
In fact,  we can prove that if $M$ is $L^{\infty}$-Liouville  then $N$ is $L^{\infty}$-Liouville.  Recalling that $M$ is $L^{\infty}$-Liouville if every bounded harmonic function $u\colon M\to\mathbb{R}$ is constant. Just lift to $M$ a harmonic bounded function $u\in C^{\infty}(N)$. The lifting $\widetilde{u}\in C^{\infty}(M)$ is bounded and harmonic thus it is constant implying that $u$ is also constant.

\subsection{Proof of Theorem \ref{Theorem 2.2}, item iii. }We start with two definitions.
\begin{definition} Let $M$ be a complete Riemannian manifold and  $v:M \rightarrow \mathbb{R}$ be a continuous function. We say that $v$ is an
exhaustion function if  all
the level sets $B^{v}_{r} = \{ x \in M; v(x) < r \} $ are
pre-compact.
\end{definition}If the exhaustion function $v$ is smooth, $C^{\infty}(M)$, then the level sets $B^{v}_{r}$ are smooth hypersurfaces for almost all $r\in v(M)\subset \mathbb{R}$.
\begin{definition} The flux of the function $v$ through a smooth
oriented hypersurface $\Gamma$ is defined by $\flux_{\Gamma} v =\displaystyle\int_{\Gamma}\langle \grad v,\,\nu\rangle d\sigma$
where $\nu$ is the outward unit  vector field normal to $\Gamma$.
\end{definition} The following theorem proved by Grigor'yan \cite[Thm. 7.6]{Grygor'yan} is fundamental in the proof of Item iii.
\begin{theorem}[Grigor'yan] A manifold $M$ is parabolic if and only if there exists a smooth exhaustion
$v$ on $M$ such that
$$\displaystyle\int_{1}^{\infty}\frac{dr}{\flux_
{\partial B^{v}_{r}}v}=\infty .$$
\end{theorem}Let $\pi \colon M\to N$ be a Riemannian manifold with compact minimal fibers $F$. Let $v\colon M\to N$ be an exhaustion function and $B^{v}_{r}$, $r>0$ its  level sets. It is clear that the lifting  $\widetilde{v}$ of $v$ is an exhaustion function of $M$ since the fibers are compact. Moreover, the level sets $B_{r}^{\widetilde{v}}$ of $\widetilde{v}$ is exactly the set $\widetilde{B^{v}_{r}}=\pi^{-1}(B^{v}_{r})=\{F_{p}=\pi^{-1}(p),\, p \in B^{v}_{r}\}$.  Let $\nu$ be the outward unit  vector field normal to $\partial B^{v}_{r}$. The lifting  $\widetilde{\nu}$ of $\nu$ is the the outward unit  vector field normal to $\partial \widetilde{B^{v}_{r}}$. Therefore, $$ \langle \grad_{M}\widetilde{v}, \widetilde{\nu}\rangle(q)=\langle  \grad_{N}v, \nu \rangle (p), \,\,\forall p\in \partial B^{v}_{r}\,\, {\rm and}\,\, \forall q\in F_{p}$$ Hence,
$$\flux_{\partial B_{r}^{\widetilde{v}} }=\int_{\partial B_{r}^{\widetilde{v}}}\langle \grad_{M}\widetilde{v}, \widetilde{\nu}\rangle\,\widetilde{d\sigma}= \int_{F_{p}}\int_{ \partial B_{r}^{v}}\langle  \grad_{N}v, \nu \rangle d\sigma(p)\,dF_{p}={\rm vol}(F_{p})\cdot \flux_{\partial B_{r}^{v}}$$
Thus,
$$\int_{1}^{\infty}\frac{dr}{\flux_{\partial B_{r}^{\widetilde{v}} }}=vol(F_{p})\cdot \int_{1}^{\infty}\frac{dr}{\flux_{\partial B_{r}^{v}}}=\infty$$ This proves that $M$ is parabolic. Observe that we used  the fact that in a Riemannian submersion with compact minimal fibers, the volume of the fibers is constant, see \cite{BMP} for more details.

\subsection{\bf Proof of Theorem \ref{Theorem 2.2}, item iv.} The proof of item iv. is an application of a recent result due to L. Mari and D. Valtorta \cite{MV} where they prove the equivalence between  the Khas'minskii criteria and stochastic completeness. We can summarize a simplified version of their result as follows.
\begin{theorem}[Khas'minskii-Mari-Valtorta]An open Riemannian manifold $M$ is stochastically complete if and only there exists a smooth exhaustion function, (called Khas'minskii function), $\gamma\colon M\to \mathbb{R}$ satisfying $\triangle \gamma \leq \lambda \, \gamma$ for some/all $\lambda >0$.
\end{theorem} By hypothesis we have a Riemannian submersion $\pi \colon M\to N$ with compact minimal fibers and the base space $N$ is stochastically complete.  Accordingly to Khas'minskii-Mari-Valtorta's Theorem there is a Khas'minskii function $\gamma$. It is straightforward to show that the lifting $\widetilde{\gamma}$ is Khas'minskii function in $M$. This shows that $M$ is stochastically complete. It should be observed that in \cite{bessa-piccione} the authors proved item iv. with  mean curvature vector of the fibers with controlled growth.

\subsection{\bf Proof of Theorem \ref{thmFellerSub}}
The idea here is to explore the relation between the minimal
solutions of the
 Dirichlet problem in $M \setminus \Omega$ and $\widetilde{M} \setminus \widetilde{\Omega}$ via an exhaustion argument which allow us to
 conclude the validity of Feller property for $M$ through the validity of the Feller property
 for $\widetilde{M}$ and reciprocally.

 Let $\{ \Omega_{n}\}_{n=1}^{\infty}$ an exhaustion of $N$ by compact sets with
 smooth boundaries. Take $\Omega \subset \Omega_{1}$ a fixed open set with smooth boundary and $\lambda >0$.
Denoting by  $\widetilde{\Omega} = \pi^{-1}(\Omega)$ and letting
$\widetilde{\Omega_{n}} = \pi^{-1}(\Omega_{n})$ we obtain an
exhaustion of $\widetilde{M}$ by compact sets with smooth boundaries.
For each $n\geq 1$, consider $h_{n}$ to be  the minimal solution of the Dirichlet problem
\begin{equation}
\left\{
\begin{array}{lll}
\triangle_{N} h_{n}  =  \lambda h, & \mbox{on} & \Omega_{n} \setminus \Omega \\
h_{n}  =  1, & \mbox{on} & \partial \Omega \\
h  =  0,  & \mbox{on} & \Omega_{n} \\
\end{array}
\right.
\end{equation}
\\
 Let
$\widetilde{h}_{n} = h_{n} \circ \pi$ be the lifting of $h_{n}$. Since the fibers are minimal, the $\widetilde{h}_{n}$'s satisfy the Dirichlet problem
\begin{equation}
%$$
\left\{
\begin{array}{lll}
\triangle_{M} \widetilde{h}_{n}  =  \lambda \widetilde{h}, & \mbox{on} & \widetilde{\Omega}_{n} \setminus \widetilde{\Omega} \\
\widetilde{h}_{n}  =  1, & \mbox{on} & \partial \widetilde{\Omega} \\
\widetilde{h}  =  0,  & \mbox{on} & \widetilde{\Omega}_{n} \\
\end{array}
\right.
%$$
\end{equation}
In fact, as the fibers $\mathcal{F}_{p}$ are minimal then we have by
\cite{BMP} that $\Div_{M}(\widetilde{X}) = \Div_{N}(X)$ where
$\widetilde{X}$ and $X$ are $\pi$-related. In particular
\begin{eqnarray*}
\triangle_{M}\widetilde{h}_{n}(\widetilde{x}) & = &
\Div_{M}(\grad_{M}\widetilde{h}_{n})(\widetilde{x}) \\
 & = & \Div_{N}(\grad_{N}h_{n})(\pi (\widetilde{x})) \\
 & = & \triangle_{N}h_{n}(\pi (\widetilde{x})) \\
 & = & \lambda h_{n}(\pi (\widetilde{x})) \\
 & = & \lambda \widetilde{h}_{n}(\widetilde{x}) \\
\end{eqnarray*}
From $\pi(\partial \pi^{-1} (\Omega)) \subset
\partial \Omega$, $\pi(\partial \pi^{-1}(\Omega_{n})) \subset
\partial {\Omega}_{n}$
we conclude that $\widetilde{h}_{n} = 1$ in $\partial
\widetilde{\Omega}$ and $\widetilde{h}_{n} = 0$ in $\partial
\widetilde{\Omega}_{n}$.  Applying  Theorem \eqref{supersol compar} we can see that for each $n\geq 1$, the functions $\widetilde{h}_{n}$ is the minimal solution for the Dirichlet problem (3.8).

Suppose that $M$ is Feller. Then we must show that given a
$\varepsilon > 0$ there exists an compact $K \subset N$ such that
$h(x) < \varepsilon$ for all $x \in N \setminus K$. Since $M$ is Feller
then given an $\varepsilon > 0$ there exists an $\widetilde{K}
\subset M$ such that $\widetilde{h}(\widetilde{x}) < \varepsilon,
\forall \widetilde{x} \in M \setminus \widetilde{K}$. Set $K =
\pi (\widetilde{K})$ and $\widetilde{K}_{0} = \pi^{-1}(K)$. We then
have that $\widetilde{K} \subset \widetilde{K}_{0}$ hence $M
\setminus \widetilde{K}_{0} \subset M \setminus \widetilde{K}$, this
means that $\widetilde{h}(\widetilde{x}) < \varepsilon, \forall
\widetilde{x} \in M \setminus \widetilde{K}_{0}$. But $\widetilde{h}
= h \circ \pi$ hence, for all $x \in N \setminus K$ we have
$$
h(x)  =  h(\pi (\widetilde{x})) = \widetilde{h}(\widetilde{x}) <  \varepsilon \\
$$
By Theorem \ref{feller equiv} we obtain that $N$ is Feller. \\
Now, suppose that $N$ is Feller, then $$\displaystyle\lim_{x \to
\infty}h(x) = 0$$ \\
 Since the fiber $\mathcal{F}_{p}$ is compact
then if $\widetilde{x} \to \infty $ on $M$ then $x \to \infty$ on
$N$ hence
$$
\displaystyle\lim_{\widetilde{x} \to
\infty}\widetilde{h}(\widetilde{x}) =
\displaystyle\lim_{\widetilde{x} \to \infty}h(\pi (\widetilde{x})) =
\displaystyle\lim_{x \to \infty} h(x) = 0
$$
that is, $M$ is Feller.

\vspace{3mm}

\noindent \textbf{ Acknowledgments:} The authors want to express their gratitude to G. Pacelli Bessa for their comments and suggestions along the preparation on this paper.

\end{document}